\documentclass[12pt]{article}

\def\tB{\tilde B}
\def\tC{\tilde C}

\def\dst{\displaystyle}
\def\hq{{\mathbf{q}}}
\def\HG#1#2#3#4#5{ \;{}_{#1} F_{#2}\left(
  {{{#3}_1,\dots,{#3}_{#1}}\atop{{#4}_1,\dots,{#4}_{#2}}}
  \,;\, #5
  \right) }

\def\qHG#1#2#3#4#5{ \;{}_{#1} \phi_{#2}\left(
  {{{#3}_1,\dots,{#3}_{#1}}\atop {{#4}_1,\dots,{#4}_{#2}}}
  \,;\, q,\,#5
 \right) }

\def\rem#1{} 

\begin{document}

\begin{center}
{\LARGE\bf
Markov's Transformation of Series \\[1ex]
and the~WZ Method
}
\\[3ex]
{\large Margo Kondratieva
and Sergey Sadov 
}
\\[3ex]
{\it
Department of Mathematics and Statistics,
Memorial University of~Newfoundland,
St.\ John's NL, A1C 5S7, Canada.
}
\\[1ex]
{E-mail: sergey@math.mun.ca}
\end{center}

\bigskip
\noindent
\hrule

\bigskip
\noindent
{\bf Abstract}\\[2ex]
In a well forgotten memoir of 1890, Andrei Markov devised a convergence
acceleration technique based on a series transformation which is very
similar to what is now known as the Wilf-Zeilberger (WZ) method.
We review Markov's work, put it in the context of modern computer-aided
WZ machinery, and speculate about possible reasons of the memoir
being shelved for so long.

\bigskip\noindent
{\it Keywords:}\
Wilf-Zeilberger method; A.A.\ Markov, Sr.; series
transformation; convergence acceleration; hypergeometric series;
basic hypergeometric series;
double series; discrete Green formula; Ap\'ery constant.

\bigskip\noindent
\hrule

\vspace{1cm}

\section{Introduction}
By this publication we aim to resurface the memoir \cite{M} by the
Russian mathematician Andrei Andreevich Markov (1856--1922), who is best
known 
as the inventor of Markov's chains in probability theory.
However, by the time Markov began his studies in probability, he was a
distinguished analyst 
and a member of the (Russian) Emperor's
Academy of Sciences.

Why would the old paper be worth attention of today's mathematical
community? All of the sudden, it appears very relevant in the
context of a powerful technique of series transformation known as
the Wilf-Zeilberger (WZ) method, and just as relevant in the
context of recent sport about faster and faster evaluation
of the constant
$$
\zeta(3)=1+\frac{1}{2^3}+\frac{1}{3^3}+\frac{1}{4^3}+\dots,
$$
called the Ap\'ery constant after R.~Ap\'ery proved its irrationality
in 1978 \cite{A}. (Not too far away is an actively pursued
challenge --- irrationality of further odd zeta values, cf. \cite{Zu} and
references therein.)

To appreciate the following results, try (if you never did) to obtain 7 correct
decimals of $\zeta(3)$ with a non-programmable 
calculator!

The formula
\begin{equation}
\label{ap}
\zeta(3)=\frac{5}{2}\sum\limits_{n=1}^{\infty}\frac{(-1)^{n-1}}{{{2n}\choose{n}}\,n^3}
\end{equation}
is often attributed to Ap\'ery, but it wasn't him who first discovered it.
The review \cite{Sri} points out the result \cite{H} reported in 1953, and here is formula (14) from Markov's memoir 
\begin{equation}
\label{ap2}
\kern-16pt
\sum_{n=0}^{\infty}\frac{1}{(a+n)^3}=\frac{1}{4}
\sum_{n=0}^{\infty} \frac{(-1)^n n!^6}{(2n+1)!}\;\frac{5(n+1)^2 +6(a-1)(n+1)
+2(a-1)^2}{[a(a+1)\dots(a+n)]^4},
\end{equation}
which is a generalization of (\ref{ap}).
The series (\ref{ap}), (\ref{ap2}) converge at the geometric rate with ratio $1/4$.
A series convergent at the geometric rate with ratio $1/27$,
\begin{equation}
\label{ap3}
\zeta(3)=\frac{1}{4}
\sum_{n=1}^{\infty} (-1)^{n-1}\, \frac{56 n^2-32 n+5}{(2n-1)^2\, n^3}\;
\frac{(n!)^3}{(3n)!},
\end{equation}
is "automatically" derived in \cite{Am}, together
with formula (\ref{ap}), using the WZ method. Interestingly,
Markov has an equivalent of (\ref{ap3}) on page 9 of his memoir.

Note for reference that \cite{Am} contains an even faster convergent representation for $\zeta(3)$ with ratio $2^2/4^4=1/64$.
A series of a non-hypergeometric type, convergent at the geometric rate with
ratio $\;e^{-2\pi}\approx 1/535\;$ is essentially due to Ramanujan
\cite[ p.30, (59)]{BBC}. And the largest number of
decimals in $\zeta(3)$, currently $520,000$, to our knowledge, was obtained
by means of the nice formula derived in \cite{AZ}
$$
\zeta(3)=\sum_{n=0}^{\infty} (-1)^n\, \frac{n!^{10}\,(205 n^2+250 n+77)}{
64\,(2n+1)!^5}.
$$
The ratio of convergence here is $2^{-10}$.

These highly nontrivial results have been obtained by the same method,
which is deceitfully simple in an abstract form. It can be
viewed either as a generalization of the elementary telescoping trick:
$$
\frac{1}{1\cdot 2}+\frac{1}{2\cdot 3}+\frac{1}{3\cdot 4}+\dots
=\left(\frac{1}{1}-\frac{1}{2}\right)+
\left(\frac{1}{2}-\frac{1}{3}\right)+
\left(\frac{1}{3}-\frac{1}{4}\right)+\dots\;=1,
$$
or as a finite-difference analog of Green's formula for
circulation of a vortex-free vector field.

One may believe in existence of interesting applications
of the discrete Green formula to series transformations,
but it isn't easy to bring forth a convincing example.
Markov has demonstrated prolificacy of that approach
in about a dozen of striking identities.
The subtlety that makes it work is a proper choice of certain
auxiliary factors unknown in advance.
(One may think of them as integrating factors).

In Sect.~2 we outline the memoir's scope and review Markov's method
(or rather its visible side).
Sect.~3 goes into details of one of Markov's examples.
We'll show close parallels between its treatment in the memoir and
by the modern computerized WZ.

That said, one should not get an impression that Markov knew the entire WZ
theory hundred years earlier. The most apparent omission in
the memoir, as well as in the later textbook \cite{MKR}, is
{\it scope}\ of the method; the related {\it how to}\
(construct such examples) and {\it what else}\ questions remain unanswered.
The creators of the modern technique put a great effort into
clarification, generalization, and algorithmization
(see \cite{PWZ,WZ90,WZ92} and
perhaps the most consonant to this context \cite{Ze}).
Also, Markov was concerned only about
convergence acceleration, while the WZ pretends to {\it
certify}, in a well-defined sense, nearly all "concrete mathematics".

In the memoir, we don't see a slightest hint to anything resembling
Gosper's algorithm (for integrating, if possible, linear difference
equations with polynomial coefficients) --- a crucial subroutine of
Zeilberger's algorithm, which, in turn, is an inborn ingredient of the WZ
method.

It would be unfair to criticize Markov for not inventing all these things.
Unfortunate --- and hard to explain --- is the fact that no one of Markov's
contemporaries picked up his technique. We speculate about
possible reasons in Sect.~4.

One of us came across the textbook \cite{MKR} in 1995 while studying
convergence acceleration methods for purposes of an applied project
\cite{Sa}. It is how the memoir \cite{M} was revealed; it is
cited in \cite{Sa}.
Unfortunately, we were not aware of the WZ method up until April 2002 and
it took more than a year for us to set on writing a detailed presentation
of Markov's work after the first published announcement \cite{KS}.
\footnote{In February 2003 Alexandru Lupas 
independently suggested  at an Internet discussion board
({\tt http://groups.google.com/groups?q=WZ-Theory})
that traces of the
WZ could be found in \cite{F,MKR} (source: \cite{Lu}).}

\section{A review of Markov's  memoir}
We begin with a translation of Section 1 of the memoir \cite{M}.

\begin{quote}
"Recall at first the proposition, which is easily derived by considering
 a double sum:

\medskip{\it
If two functions $\,U_{x,z}$ and $V_{x,z}$ of independent variables $x$ and $z$
are bound by the condition
\begin{equation}
\label{disdiv}
 U_{x,z}-U_{x+1,z}=V_{x,z}-V_{x,z+1},
\end{equation}
then
\begin{equation}
\label{disgreen}
\begin{array}{l}\dst
 U_{0,0}+U_{0,1}+\dots +U_{0,j-1} \; - \;U_{i,0}- \;U_{i,1}-\dots -\;U_{i,j-1}
 \\[2ex]
 \qquad\dst =\;
 V_{0,0}+V_{1,0}+\dots +V_{i-1,0}\;-\;V_{0,j}-V_{1,j}-\dots -V_{i-1,j},
\end{array}
\end{equation}
$i$ and $j$ being arbitrary positive integers.
}

\medskip
In all the cases occurring in this memoir, the series with terms
$$
\begin{array}{l}
 U_{0,0}, U_{0,1},\dots,U_{0,j},\dots
 \\[1.2ex]
 V_{0,0}, V_{1,0},\dots,U_{i,0},\dots
\end{array}
$$
are convergent and the sums
$$
 U_{i,0}+U_{i,1}+ \dots +U_{i,j-1},\dots
 \qquad{\rm and}\qquad
 V_{0,j}, V_{1,j}+ \dots +U_{i-1,j},\dots
$$
tend to zero as $i$ and $j$ increase indefinitely.

That stated, the formula (\ref{disgreen}) will give
\begin{equation}
\label{martran}
 U_{0,0}+U_{0,1}+\dots +U_{0,j} +\dots\;=\;
 V_{0,0}+V_{1,0}+\dots +V_{i,0}+\dots ."
\end{equation}
\end{quote}

Markov works with hypergeometric and basic hypergeometric
series in his memoir, although he avoids calling them so.
To make formulae more concise and comprehensible,
let us recall appropriate definitions and notation, cf.\ \cite{Ba,GR}.

The {\it rising factorial} is defined as
$$
(a)_n=a(a+1)\cdot\dots\cdot (a+n-1)=\frac{\Gamma(a+n)}{\Gamma(a)},
\qquad n\geq 0.
$$
In particular, $(1)_n=n!$.
Denote for brevity
$$
 (a_1,\dots,a_r)_n=\prod_{j=1}^r (a_j)_n.
$$
A {\it hypergeometric {\rm (HG)} term}  is an expression of the form
$$
\frac{(a_1,\dots,a_r)_n}{(b_1,\dots,b_s)_n} \;z^n,
$$
and a hypergeometric series is a series of the form
$$
 \HG{r}{s}{a}{b}{z}=\sum_{n=0}^{\infty} \frac{(a_1,\dots,a_r)_n}{(b_1,\dots,b_s,1)_n}\; z^n
$$
If $z=1$, it is common to omit the argument $z$ .

Basic hypergeometric {\rm (BHG)} terms and series contain an additional parameter
$q$,  
called the {\it base}.
The {\it $q$-rising factorial}\ is the product
$$
(a;q)_n =(1-a)(1-qa)\cdot\dots\cdot(1-q^{n-1}a),
$$
The expression $(q;q)_n$ is called the {\it $q$-factorial} of $n$.
The product of several $q$-rising factorials is abbreviated as
$$
(a_1,\dots,a_r;\,q)_n=\prod_{j=1}^r (a_j;q)_n.
$$
A {\it basic hypergeometric term}\ is an expression of the form
$$
\frac{(a_1,\dots,a_r;\,q)_n}{(b_1,\dots,b_s;\,q)_n}\; z^n,
$$
and a basic hypergeometric series is a series of the form
$$
 \qHG{r}{s}{a}{b}{z}=\sum_{n=0}^{\infty} \frac{(a_1,\dots,a_r;\,q)_n}{(b_1,\dots,b_s,1;\,q)_n}
  \;z^n \left((-1)^n q^{n(n-1)/2}\right)^{1+s-r}.
$$
The base $q$ is usually omitted in the notation, unless
BHG series with different bases are discussed in the same context.
Also, as for HG series, the argument $z$ is often omitted in the special
case $z=1$.

The ordinary hypergeometry is a limiting case of the basic one:
\begin{equation}
\label{limBHG}
\lim_{q\to 1}\,\frac{(q^a;q)_n}{(q^b;q)_n}=\frac{(a)_n}{(b)_n}.
\end{equation}
In the basic case, Markov assumes $\,|q|>1$,
while the modern convention strongly prefers $\,|q|<1$.
For this reason we re-denote Markov's $q\,$ to $\;\hq\;$
and adopt the base $\;q=\hq^{-1}$.
Thus in the sequel $|q|<1$ and $|\hq|>1$.

\bigskip\noindent
{\bf Structure of the functions $U_{x,z}$ and $V_{x,z}$ in Markov's examples}
\\[2ex]
All the examples, in their most general form, deal with convergence-accelerat\-ing
transformations of  hypergeometric or basic hypergeometric series. Every time
we have a HG or BHG term
$F_{x,z}$ such that the series $\sum_z F_{0,z}$ is to be summed.
Dependence of $F_{x,z}$ on $x$ is characterized by the multiplicative pattern
$$
\frac{(a)_z}{(b)_{x+z}};
$$
in the basic case an additional factor $q^{f(x,z)}$ is present,
where $f$ is a quadratic polynomial in $x$.
In all cases, $F_{x,z}$ and the sums over $z$ have limit $0$ as $x\to\infty$.

The function
$U_{x,z}$ has one of the following three forms:
\begin{equation}
\label{u1}
 U_{x,z} = F_{x,z} A_{x}\qquad\qquad \mbox{\rm in \S\,2,3}
\end{equation}
or
\begin{equation}
\label{u2}
 U_{x,z} = F_{x,z} (A_{x}+B_{x}z) \qquad\quad \mbox{\rm in \S\,4,8}
\end{equation}
or
\begin{equation}
\label{u3}
 U_{x,z} = F_{x,z} (A_{x}+\tB_{x}z+\tC_{x}z^2) \qquad \mbox{\rm in \S\,9}.
\end{equation}
In the first case $F_{x,z}$ is a BHG
term, and in the latter two cases $F_{x,z}$ is a HG term.
In addition, 
$B_0=0$ and $\tB_0=\tC_0=0$.
Trying to present Markov's patterns in a unified form,
we use symbols with tildes where our notation is not identical
to that in \cite{M}.

The function $V_{x,z}$ is sought in the form
\begin{equation}
\label{vm}
 V_{x,z}=F_{x,z} M_{x,z},
\end{equation}
where in the case (\ref{u1})
$$
M_{x,z}=B_x+C_x q^{z};
\eqno(\ref{u1}')
$$
in the case (\ref{u2})
$$
M_{x,z}=C_x+ D_x z+ F_x z^2,
\eqno(\ref{u2}')
$$
and finally in the case (\ref{u3})
$$
M_{x,z}=F_{x}+\tilde G_{x}z+\tilde H_{x}z^2.
\eqno(\ref{u3}')
$$

\subsection*{\bf List of the series dealt with in the Memoir}

\hangindent 20pt\noindent
\S\,2:
$\qquad\dst{}_2\phi_1\left({{a,1}\atop{b}}; t\right)\,
=\,1\,+\frac{1-a}{1-b}\,t\,+\frac{(1-a)(1-aq)}{(1-b)(1-bq)}\,t^2\,+\,\cdots
$.

\hangindent 20pt\noindent
\S\,3:$\qquad\dst{}_3\phi_2\left({{a,b,1} \atop{c,d}};
\dst\frac{cd}{abq}\right)\,$ and the limiting (Schellbach's) case
$\dst\;\,{}_3 F_2\left({{a,b,1}\atop {c,d}};\right)\,$.

\hangindent 20pt\noindent
\S\,4:$\qquad\dst{}_4 F_3\left({{a,a+h,a-h,1}\atop{b,b+h,b-h}};\right)\,
=\,1\,+\,\frac{a}{b}\,\cdot\,\frac{a^2-h^2}{b^2-h^2}$ \\[2ex]
\hspace*{5cm}
$\dst\,+\,
\frac{a}{b}\,\cdot\,\frac{a+1}{b+1}\,\cdot\,\frac{a^2-h^2}{b^2-h^2}\,
\cdot\,\frac{(a+1)^2-h^2}{(b+1)^2-h^2}\,+\,\cdots
$.

\hangindent 20pt\noindent
\S\,5: $\;$ Special case of \S\,4:
$\;a=1,\;b=2$; then further specialization $\,h=0\,$
yields the series defining $\zeta(3)$. Formula (\ref{ap3}) is
found in this \S.

\hangindent 20pt\noindent
\S\,6: $\;$ Special case of \S\,4: $h=0,\;b=a+1$, yielding the
Hurwitz zeta series $\zeta(3,a)\,$.

\hangindent 20pt\noindent
\S\,7: $\;$ Special case of \S\,4: Kummer's sum
$$
{}_4 F_3\left({{\frac{9}{2},\,\frac{9}{2},\,\frac{9}{2},\,1}
\atop{5,\,5,\,5}};\right)
\,=\sum_{n=0}^{\infty}
\left(\frac{(2n+1)!!}{(2n)!!}\right)^3.
$$

\hangindent 20pt\noindent
\S\,8: $\;$ The series
$\dst\;{}_3 F_2\left({{a,b,1}\atop {c,d}};\,-1\right)\;$
with $\,c-a=d-b$, and particular cases.

\hangindent 20pt\noindent
\S\,9: $\;$ A well-poised \cite{Ba,GR} series
$\dst\;{}_4 F_3\left({{a,a,a,1}\atop{b,b,b}}; \,-1\right)\,$.
Among considered particular cases there are Stirling's series
$\sum_1^\infty (-1)^n n^{-k}$, $\;k=2,3$.

\noindent
Three more formulae, (\ref{ap2}) being the simplest, are contained
in the last \S\,10. Details of the transformations, in particular,
the form of the functions $U$, $V$, are not provided.

\subsection*{\bf Method for obtaining the transformations}
Once the parametric form of the functions $U_{x,z}$ and
$V_{x,z}$ is set, it remains to choose the undetermined constants
in order to satisfy Eq.\ (\ref{disdiv}). The main question
is why exactly that many parameters are needed in the particular
situation. We suppose that Markov simply used a trial and error approach,
starting with minimal number of parameters and extending the family
of parameters until a solution was found. But there may exist
a clever reasoning of which we are not aware.
\footnote{Written before we had a chance to study \cite{MZ}.}

\section{Example: Transformation of a $\,_3\phi_2\,$ series}
Consider the series
\begin{equation}
\label{qsh}
_3\phi_2\left({{a,b,1}\atop{c,d}};\, q,\,t\right)=
\sum_{z=0}^\infty \frac{(a;q)_z (b;q)_z}{(c;q)_z (d;q)_z}\,t^z,
\qquad t=\frac{cd}{abq},
\quad |t|<1.
\end{equation}
from \S\,3 of Markov's memoir.
We made a 
random choice between this example and the one in \S\,2, but
we deliberately chose an example that falls under the case (\ref{u1}),
where the auxiliary factor is $A_x$ is $z$-independent.
Our intention is to compare Markov's procedure with the WZ one.

Markov writes the series in the form (as agreed, we rename his $q$ to
$\hq=q^{-1}$)
\begin{equation}
\label{qshm}
1+\frac{(r-1)(r'-1)}{(s-1)(s'-1)}\,\hq\,+\,
\frac{(r-1)(r\hq-1)(r'-1)(r'\hq-1)}{(s-1)(s\hq-1)(s'-1)(s'\hq-1)}\,\hq^2\,
+\dots.
\end{equation}
In view of the relation
$$\;
(u-1)(u\hq-1)\dots (u\hq^{n-1}-1)=(u^{-1};q)_n \;u^n\, q^{n(n-1)/2},
$$
the series (\ref{qsh}) and (\ref{qshm}) are equivalent, if their
parameters are related as follows
$$
 a=\frac 1{r}, \quad\; b=\frac 1{r'},\quad\;c=\frac 1{s}, \quad\; d=\frac 1{s'}, \quad\; t=\frac{cd}{abq}=\frac{rr'\hq}{ss'}.
$$
First, let us derive Markov's result (see (\ref{v0}) below) in our notation by hand.

Set $A_0=1$ in (\ref{u1}). Markov takes an extension
$F_{x,z}$ of the term $\;F_{0,z}=U_{0,z}$ $= t^{z}(a,b;q)_z /(c,d;q)_z \;$
in the form
\begin{equation}
\label{F}
 F_{x,z}= \frac{(a;q)_z (b;q)_z \,t^{z}}{(c;q)_{x+z} (d;q)_{x+z}} \;(cd q^{2z})^x \,q^{x(x-1)}.
\end{equation}
Such a pattern doesn't appear obvious when using the modern
form with $|q|<1$, but it is naturally suggested by the original form
(\ref{qshm}): replace the two products of $z$ factors in the
denominator by the products of $(x+z)$ factors.  In particular,
$$
\;F_{x,0}=\,q^{x(x-1)}\;\frac{(cd)^{x}}{(c,d;q)_{x}}.
$$

With (\ref{F}) and (\ref{u1}$'$), the condition (\ref{disdiv}) becomes
$$
\begin{array}{l}\dst
A_{x}\,\frac{(a,b;q)_z\,t^z}{(c,d;q)_{x+z}}\,(cdq^{2z})^x\,q^{x(x-1)}\,
-\,A_{x+1}\,\frac{(a,b;q)_z\,t^z}{(c,d;q)_{x+z+1}}\,(cdq^{2z})^{x+1}\,q^{(x+1)x}
\\[3ex]\dst
=
\left[(B_x+C_{x}q^z)\,\frac{(a,b;q)_z\,t^z}{(c,d;q)_{x+z}}\,-\,
(B_{x}+C_{x}q^{z+1})\,\frac{(a,b;q)_{z+1}\,t^{z+1}}{(c,d;q)_{x+z+1}}\,q^{2x}\,\right]\,
\\[2.5ex]\dst \qquad
\times
\,(cdq^{2z+x-1})^x.
\end{array}
$$
Taking out the common factor
$\,(a,b;q)_z/(c,d;q)_{x+z+1}\,t^z\,(cdq^{2z+x-1})^x$, we
obtain an equation of degree 3 in $q^z$. To satisfy condition
(\ref{disdiv}), all the coefficients of that equation must vanish,
that is
{
\begin{eqnarray}
 \dst A_x &=&\dst B_x(1-tq^{2x}),\label{c0}\\[1.8ex]
 \dst -A_x(c+d)q^x &=&\dst C_x-B_x(c+d)q^x+B_x(a+b)q^{2x}t-C_x q^{2x+1}t,\label{c1}\\[1.8ex]
 \kern-14pt\dst (A_x -A_{x+1}) cd\,q^{2x}&=&\dst B_x
 (cd-ab t)\,q^{2x}+C_x
 \left((a+b)q^{2x+1}t-(c+d)q^x\right)
 ,\label{c2}\\[1.8ex]
 \dst 0 &=&\dst C_x(cd\,q^{2x}-ab\,q^{2x+1}\,t)\label{c3}.
\end{eqnarray}
}
Eq.\ (\ref{c3}) holds automatically. Equations (\ref{c0}),
(\ref{c1}) imply
\begin{equation}
\label{abc}
\begin{array}{l}\dst
\frac{B_{x}}{A_x}=\left(1-tq^{2x}\right)^{-1},
\\[2ex]
\dst
\frac{C_x}{A_x}=\frac{tq^{2x}[(c+d)q^{x}-(a+b)]}{
(1-tq^{2x})(1-tq^{2x+1})}
\end{array}
\end{equation}
and therefore, by (\ref{u1}')
\begin{equation}
\label{ma}
 \frac{M_{x,0}}{A_x}=\frac{B_x+q^0 C_x}{A_x}
 =\frac{1-tq^{2x}(a+b+q)+tq^{3x}(c+d)}{
 (1-tq^{2x})(1-tq^{2x+1})}.
\end{equation}
Substitution of (\ref{abc}) to (\ref{c2}) yields the recurrence
for $A_x$
\begin{equation}
\label{aa}
 \frac{A_{x+1}}{A_x}=\frac{(1-\frac{c}{a} q^x)(1-\frac{c}{b}q^x)
 (1-\frac{d}{a} q^x)(1-\frac{d}{b}q^x)}{q\,(1-tq^{2x})(1-tq^{2x+1})}.
\end{equation}
Since $A_0=1$, we obtain
$$
A_x=\frac{(\frac{c}{a}, \frac{c}{b}, \frac{d}{a}, \frac{d}{b};\, q)_x}{q^x\,(t;q)_{2x}}.
$$
Finally, we find the general term of the transformed series in the
r.h.s.\ of (\ref{martran})
\begin{equation}
\label{v0}
\begin{array}{rcl}
 \kern-14pt
 V_{x,0}&=&\dst\frac{M_{x,0}}{A_x} A_x F_{x,0}
 \\[2.5ex]
 &=&\dst
 \frac{(\frac{c}{a}, \frac{c}{b}, \frac{d}{a}, \frac{d}{b};\, q)_x}{
 (c,d;\,q)_x} \; (cd)^x\,q^{x(x-2)} \;
 \frac{1-tq^{2x}(a+b+q)+tq^{3x}(c+d)}{(t;\,q)_{2x+2}}.
 \end{array}
\end{equation}
For $x=0$ or $1$, the terms are consistent with those in formula
(7) in \cite{M}, where further terms are not written out, while
they are not easy to guess.

Following Markov, we proceed to consider the limiting case $q\to 1$.
Re-denote $a,b,c,d,t$\ respectively to $q^a, q^b, q^c, q^d, q^t$.
The relation $tq=(cd)/(ab)$ is replaced by the following:
$$
 t+1=c+d-a-b.
$$
Then, by (\ref{limBHG}),
$$
\lim_{q\to 1} \frac{(q^{c-a},q^{c-b},q^{d-a},q^{d-b};\,q)_x}{(q^{c},q^{d};\,q)_x\,(t;q)_{2x}}
\,=\,\frac{({c-a},{c-b},{d-a},{d-b})_x}{({c},{d})_x\,(t)_{2x}}.
$$
The limit of the remaining factor in (\ref{v0}) is found by applying L'Hospital Rule two times:
$$
\lim_{q\to 1} \frac{1-q^{2x+t}(q^a+q^b+q)+q^{3x+t}(q^c+q^d)}{(1-q^{2x+t})(1-q^{2x+t+1})}
\,=\,\frac{p(a,b,c,d,x)}{(2x+t)(2x+t+1)},
$$
where
$$
\begin{array}{l}
p(a,b,c,d,x)\,=\,
(2x+t+a)(2x+b+t)-(x+c-1)(x+d-1)\\[1.7ex]
\qquad
=\, (c+d-a-1+2x)(c+d-b-1+2x)-(c-1+x)(d-1+x).
\end{array}
$$
The result is Schellbach's formula
$$
{}_3 F_2\left({{a,b,1}\atop{c,d}};\right)=\,\sum_{x=0}^{\infty}
\frac{({c-a},{c-b},{d-a},{d-b})_x\;p(a,b,c,d,x)}{({c},{d})_x\,(c+d-a-b-1)_{2x+2}}.
$$
The left-hand side converges as $\sum n^{-t-1}$
(assuming that $t=c+d-a-b-1>0$).
The right-hand side converges geometrically; namely,
the term with subscript $x$ has the asymptotics
$\,4^{-x}\cdot K\,x^{-a-b+1/2}\;$
with
$$\,K=(3/2^{t+1})\sqrt{\pi}\Gamma(t)\Gamma(c)\Gamma(d)/
(\Gamma(c-a)\Gamma(c-b)\Gamma(d-a)\Gamma(d-b)).
$$

We turn now to the Wilf-Zeilberger approach, more specifically, to its computer-aided version. Speaking pragmatically,
all one needs is to type in the expression (\ref{F}) in Maple, feed it to the {\tt qEKHAD} program, and analyze the results.
The substitution $t=(cd)/(abq)$ must be made in advance in (\ref{F}).

The program produces a {\it recurrence operator} $\,\Omega(X,x)$ and a {\it certificate} $R(x,z)$. We believe that the reader
can't avoid looking into \cite{PWZ} anyway, but below we give a
self-contained account of the procedure in this case.

The recurrence operator outputted by {\tt qEKHAD} has the structure
$$
\Omega(X,x)=P(x)+Q(x)X.
$$
Here $X$ is the operator of forward shift in $x$, that is
\begin{equation}
\label{rec}
 (\Omega(X,x)F)_{x,z}=P(x)F_{x,z}+Q(x)F_{x+1,z}.
\end{equation}
The certificate $R(x,z)$ is a rational function of $q^x$, $q^z$ such that the function
$$
 G(x,z)=R(x,z)F_{x,z}
$$
satisfies the equation
$$
 (\Omega(X,x)F)_{x,z}=G(x,z)-G(x,z-1).
$$
For comparison purposes, it is more convenient to deal with {\it forward}\
\linebreak[4]
$z$-difference in the right-hand side, so we denote
$$
 \tilde G(x,z)=G(x,z-1)=\tilde R(x,z)F_{x,z},
$$
where
$$
\tilde R(x,z)=R(x,z-1)\frac{F_{x,z-1}}{F_{x,z}}.
$$
Now
\begin{equation}
\label{om1}
 (\Omega(X,x)F)_{x,z}=\tilde G(x,z+1)-\tilde G(x,z).
\end{equation}
We will actually need only values $\tilde R(x,0)$.
Taking the output of {\tt qEKHAD} and transforming it this way, we
find (with $t=(cd)/(abq)$, as before)
\begin{equation}
\label{r0}
\tilde R(x,0)=\frac{1-tq^{2x}(a+b+q)+ tq^{3x}(c+d)}{
1-t q^{2x+1}}.
\end{equation}
The values of $P(x)$ and $Q(x)$ in (\ref{rec}) produced by {\tt qEKHAD} are
\begin{equation}
\label{p0}
P(x)=1-tq^{2x},\qquad Q(x)=\frac{(1-\frac{c}{a}q^x)(1-\frac{c}{b}q^x)
(1-\frac{d}{a}q^x)(1-\frac{d}{b}q^x)}{q(1-t q^{2x+1})}.
\end{equation}
Equations (\ref{r0}), (\ref{p0}) have much in common with (\ref{ma}), (\ref{aa}), though they are not identical. Of course, the similarity is not occasional. It is explored below in detail.

If we fix $F_{x,z}$ and try to satisfy Eq.~(\ref{disdiv}) using substitutions of the form (\ref{u1}), (\ref{vm}), the following equation comes up:
\begin{equation}
\label{rec2}
 A_{x+1} F_{x+1,z}-A_x F_{x,z}\,=\, M_{x,z+1} F_{x,z+1}-M_{x,z} F_{x,z}.
\end{equation}
On the other hand, Eq.~(\ref{om1}) in expanded notation reads
\begin{equation}
\label{om2}
 Q(x) F_{x+1,z}+P(x)F_{x,z}\,=\,
 \tilde R(x,z+1) F_{x,z+1}-\tilde R(x,z) F_{x,z}.
\end{equation}
Suppose that the certificate $\tilde R(x,z)$ and the operators $\Omega(x,z)$ are known. Let us find $A_x$ and $M_{x,z}$. Introduce as yet undetermined coefficients $\Phi(x)$ such that multiplication by $\Phi(x)$ turns Eq.~(\ref{om2}) into (\ref{rec2}). Thus,
$$
 M_{x,z}=\Phi(x)\tilde R(x,z)
$$
and
$$
A_x=\Phi(x) P(x),\qquad A_{x+1}=\Phi(x) Q(x).
$$
Therefore,
$$
\frac{M_{x,z}}{A_x}=\frac{\tilde R(x,z)}{P(x)},\qquad  
\quad
\frac{A_{x+1}}{A_x}=\frac{Q(x)}{P(x)}.
$$
The right-hand sides in these equations follow from (\ref{r0}), (\ref{p0}).
The obtained equations for $A_x$ and $M_{x,z}$ are identical to (\ref{ma}) and (\ref{aa}),
from which we (following Markov) have found the terms (\ref{v0}) of the transformed series.

\section{How did Markov miss his audience?}

This section is mostly speculative. A thorough study of
Markov's works, letters, and other documents, which may
reveal circumstances of the appearance of the memoir in question and
of its abandonment, is yet to be undertaken.

Having been deeply involved in studies on continued fractions throughout
the 1880s, Markov corresponded with T.~J.~Stieltjes 
\cite{OS} and closely
watched his publications. In 1887 Stieltjes \cite{Sti} published a table
of the values of the Riemann Zeta function $\zeta(k)$ with 32 decimals
for integral values of $k$ from 2 to 70. Markov might have felt challenged
by that achievement and by Stieltjes' convergence acceleration technique.
Apparently, it was this challenge and rivalry that prompted Markov to develop
his new acceleration method.
In a brief note \cite{M89}
he gives two formulae, one of them equivalent to (\ref{ap3}),
and obtains 20 decimals of $\zeta(3)$ taking 13 terms
in his series.
Afterwards he jealously beat Stieltjes' record,
taking 22 terms and obtaining the result with 33 decimals in \cite{M}
$$
 \zeta(3)= 1,20205 69031 59594 28539 97381 61511 450.
$$
The second formula published in \cite{M89} is a $\,27^{-k}\,$-fast convergent
representation
$$
\zeta(2)=\frac{5}{3}+\sum_{k=1}^{\infty}\frac{(-1)^k (2k-1)!!^3}{(6k-1)!!}
\left(\frac{1}{4k^2}+\frac{5}{(6k+1)(6k+3)}\right).
$$

Claiming that Markov missed his audience, as it eventually turned out,
we don't mean that the series transformation he proposed remained unnoticed.
Markov himself tried to popularize it. In the textbook \cite{MKR}
there is a chapter devoted to this transformation with a number of examples,
although examples with basic hypergeometric series are not included.
References to Markov's work are found in the well-known textbooks
\cite[III.24]{B},
\cite{K}. The latter contains a section (Ch.~VIII, \S~33) on Markov's
transformation, at the beginning of which we find, among all, a reference
to Stirling's work 
[27], the starting point of Gosper's seminal
paper \cite{Gos}, which laid out the foundation of automated identity proving.
Both Stirling's and Markov's methods are treated in detail in another text
\cite{F}, which seems to be left out nowadays, perhaps undeservedly.

The evidence that T.~J.~Bromwich \cite{B} was aware of Markov's work is especially interesting,
since it was England where the research in hypergeometric and combinatorial
identities
enjoyed its most fruitful period in the first two decades of the 20th
century.  Did Rogers and MacMahon see Markov's memoir? Ramanujan might have
appreciated formula (\ref{ap}) had he noticed it in \cite{B},
but Hardy \cite[II.14]{Ha} doubts Ramanujan having seen that book.

It is perhaps even more surprising that the memoir of 1890 had been completely
forgotten in Russia. Markov's name and works were well known and highly regarded
in the Soviet Union.%
\footnote{To avoid a confusion, we are talking about A.~A.~Markov, Sr.
His son, Andrei Andreevich Markov, Jr.\ (1903--1979), was also a
prominent mathematician, a member of the USSR Academy of Sciences, and
one of the founders of Computer Science in the Soviet Union.
}
The biography \cite{Gro} contains an appendix, where Markov's works in various
directions are reviewed by experts in the respected areas. In the Analysis section
(as well as anywhere else), the convergence acceleration topic is not even
mentioned! 
We managed to find only one reference to \cite{M} in mathematical literature
of the Soviet period: a rarity textbook \cite{Rom}. It contains a section
on Markov's transformation and the exposition there, as the author indicates,
closely follows that in \cite{K}. "Markov's theorem", see below,
is also found in a widely circulated treatise \cite{Vor}; however, no exact
reference and no applications are given.

In our opinion, the latter theorem is partly to blame for
draining the key issue of the 1890 memoir. The theorem is also contained in
the cited texts \cite{K} and \cite{Rom}, and it goes back to Markov's
lecture notes \cite{MKR}. The formulation below is taken from \cite{K}.

\medskip\noindent
{\bf Theorem.} {\it
Let a convergent series $\sum_{k=0}^{\infty} z^{(k)}$ be given with each of its
terms itself expressed as a convergent series:
\begin{equation}
\label{mth}
z^{(k)}=a_0^{(k)}+a_1^{(k)}+\dots+a_n^{(k)}+\dots \qquad
(k=0,1,2,\dots).
\end{equation}
Let the individual columns $\;\dst\sum_{k=0}^{\infty} a_n^{(k)}\;$
of the array {\rm(\ref{mth})} so formed represent convergent series with sum
$s^{(n)}$, $\;n=0,1,2,\dots$, so that the remainders
$$
 r_m^{(k)}=\sum_{n=m}^{\infty}a_n^{(k)} \qquad (m\geq 0)
$$
of the series in the horizontal rows also constitute a convergent series
$$
 \sum_{k=0}^{\infty} r_m^{(k)}=R_m \qquad (m \;\;{\it fixed}\,).
$$
In order that the sums by vertical columns should form a convergent series
$\,\sum s^{(n)}$, it is necessary and sufficient that $\;\lim R_m=R\,$ should
exist; and in order that the relation
\begin{equation}\label{mt2}
\sum_{n=0}^{\infty} s^{(n)} =\sum_{k=0}^{\infty} z^{(k)}
\end{equation}
should hold as well, it is necessary and sufficient that this limit should
be $0$.
}

\medskip
Compare the introductory section of the memoir and this
theorem. The latter is as simple in essence as the former but how much
harder it is to grasp! It is positioned, in the first instance, as a {\it convergence
theorem} and the equation (\ref{mt2}) is just yet another switch-the-order
formula. The theorem {\it per se} expresses a nice and possibly useful
analytical criterion, but it completely overshadows the original point.

\medskip
This may partly explain an underestimation of Markov's work, but another
component is the strikingly different level (compared to the nearly trivial
general idea) of concrete formulae, and a lack of Markov's elaboration
on the forms of the series and the auxiliary factors.
In neither of the cited books did their authors offer their own examples!
And, since Markov didn't make any precise statements regarding
the applicability range of his transformation, we got to observe
the tendency to phase out vague and complicated applications and
emphasize the simple and well-rounded theorem. But does calculus
exist for the sake of convergence theorems?

Ch.~Hermite, the then-editor of {\it Comptes Rendus,}\ replied to
\cite{M89}:
{\it "Par quelle voie vous \^etes parvenue \`a une telle transformation,
je ne puis m\^eme de loin l'entrevoir,
et il me faut vous laisser v\^otre secret."}
\footnote{"I can't even remotely guess the way you arrived at such a transformation,
and it remains to leave your secret with you."}
\cite{O},\cite{Todd}
It might have sound as a compliment, but shouldn't it be heard by
Markov as a warning?
A nice hint would have helped Markov's readers (and himself?): an
advise to investigate specifically series of hypergeometric and basic
hypergeometric type.
Restriction to these two classes yielded the development of effective
algorithms whose traces are implicit in \cite{M} for determining
auxiliary factors (certificates, in the WZ version) and ensured a
huge success of the modern WZ method.
Availability of a software ({\tt EKHAD}, {\tt qEKHAD} --- see \cite{PWZ})
makes it tremendously helpful for everybody who deals with hypergeometric
functions, partitions, and the like.

\section*{Acknowledgments}

We are grateful to Natalia Ermolaeva for helpful information about
\linebreak[4]
A.~A.~Markov and in particular for pointing out Ref.~\cite{Rom}.
A stimulating inquiry by Victor Adamchik prompted us
to finally write the long planned text.
We appreciate constructive comments of anonymous referees
for the ISSAC'04 Conference, where the paper was originally intended to.
It is our pleasure to thank Doron Zeilberger for his interest that
(understandably) went far beyond editorial remarks and has been
expressed by now in the article \cite{MZ}.

M.K.~acknowledges a support by a grant from the
Natural Sciences and Engineering Research Council of Canada.

\newcommand{\iT}{}
\newcommand{\bF}{}

\end{document}